\documentclass[runningheads]{llncs}

\def\anonymous{no}

\usepackage{ifthen}

\usepackage{amsmath,amssymb}
\usepackage{enumitem}

\usepackage{color}
\usepackage[usenames,dvipsnames]{xcolor}

\usepackage{graphicx}
\usepackage{subcaption}
\usepackage[unicode,colorlinks=true,urlcolor=blue,citecolor=red]{hyperref}%
\graphicspath{{./figs/}}
\usepackage[capitalise, noabbrev]{cleveref}

\usepackage[mathlines]{lineno}

\usepackage{csquotes}

\usepackage{todonotes}

\spnewtheorem{observation}[theorem]{Observation}{\bfseries}{\itshape}
\spnewtheorem{ques}{Question}{\bfseries}{\itshape}
\usepackage{thmtools, thm-restate}

\renewcommand{\orcidID}[1]{\href{https://orcid.org/#1}{\includegraphics[scale=.03]{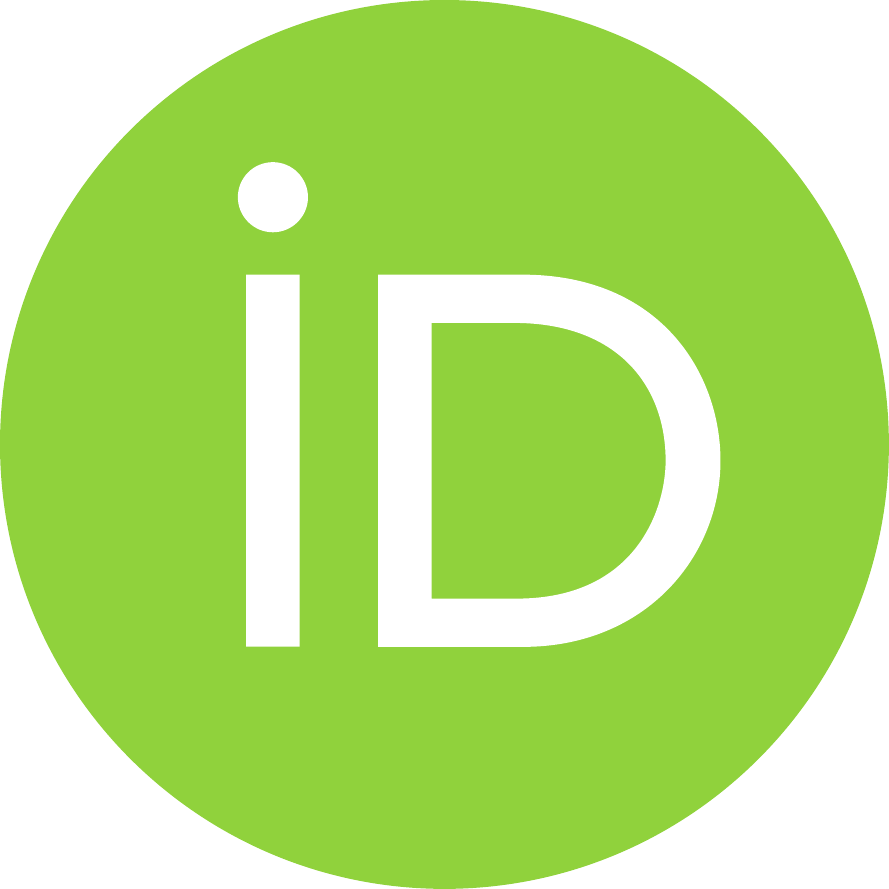}}}

\begin{document}

\ifthenelse{\equal{\anonymous}{no}}{

\title{Compatible Spanning Trees in Simple Drawings of $K_n$ \thanks{This work was initiated at the 6th DACH Workshop on Arrangements and Drawings in Stels, August 2021. We thank all participants, especially Nicolas~Grelier and Daniel~Perz, for fruitful discussions. O.A., R.P.~and A.W.~are supported by FWF grant~W1230. K.K.~is supported by the German Science Foundation (DFG) within the research training group `Facets of Complexity' (GRK 2434). W.M.~is partially supported by the German Research Foundation within the collaborative DACH project \emph{Arrangements and Drawings} as DFG Project MU 3501/3-1, and by ERC StG 757609. J.O.~is supported by ERC StG 757609. M.M.R.~is supported by the Swiss National Science Foundation within the collaborative DACH project \emph{Arrangements and Drawings} as SNSF Project 200021E-171681. (Also note that this author's full last name consists of two words and is \emph{Mallik Reddy}. However, she consistently refers to herself with the first word of her last name being abbreviated.) B.V.\ was partially supported by the Austrian Science Fund (FWF) within the collaborative DACH project \emph{Arrangements and Drawings} as FWF project \mbox{I 3340-N35}.}}
\author{Oswin~Aichholzer\inst{1}\orcidID{0000-0002-2364-0583} \and
	Kristin~Knorr\inst{2}\orcidID{0000-0003-4239-424X} \and
	Wolfgang~Mulzer\inst{2} \and 
	Nicolas~El~Maalouly\inst{3}\orcidID{0000-0002-1037-0203} \and
    Johannes~Obenaus\inst{2}\orcidID{0000-0002-0179-125X} \and
	Rosna~Paul\inst{1}\orcidID{0000-0002-2458-6427} \and
	Meghana~M.~Reddy\inst{3}\orcidID{0000-0001-9185-1246} \and
	Birgit~Vogtenhuber\inst{1}\orcidID{0000-0002-7166-4467} \and
	Alexandra~Weinberger\inst{1}\orcidID{0000-0001-8553-6661}
}
\authorrunning{O. Aichholzer et al.}

\institute{Institute of Software Technology, Graz University of Technology, Austria\\
	\email{\{oaich,ropaul,bvogt,weinberger\}@ist.tugraz.at} \and
	Institut für Informatik, Freie Universit{\"a}t Berlin, Germany\\
	\email{\{kristin.knorr,wolfgang.mulzer,johannes.obenaus\}@fu-berlin.de} \and
	 Department of Computer Science, ETH Z{\"u}rich, Switzerland\\
	\email{\{nicolas.elmaalouly,meghana.mreddy\}@inf.ethz.ch}
}

}{
\title{Compatible Spanning Trees in Simple Drawings of $K_n$}
\author{\textcolor{red}{anonymous author(s)}}
\authorrunning{anonymous}
}

\maketitle

%\linenumbers

\begin{abstract}
	For a simple drawing $D$ of the complete graph $K_n$, two (plane) subdrawings are \emph{compatible} if their union is plane. Let $\mathcal{T}_D$ be the set of all plane spanning trees on $D$ and $\mathcal{F}(\mathcal{T}_D)$ be the \emph{compatibility graph} that has a vertex for  
each element in $\mathcal{T}_D$ and two vertices are adjacent if and only if the corresponding trees are compatible. We show, on the one hand, that $\mathcal{F}(\mathcal{T}_D)$ is connected if $D$ is a cylindrical, monotone, or strongly c-monotone drawing. On the other hand, we show that the subgraph of $\mathcal{F}(\mathcal{T}_D)$ induced by stars, double stars, and twin stars is also connected.
	In all cases the diameter of 
	the corresponding compatibility graph is at most linear in $n$.
 
\end{abstract}

\section{Introduction}

A \emph{drawing} $D$ of a graph $G$ is a representation of $G$ in the Euclidean plane such that the vertices of $G$ are distinct points and the edges are Jordan arcs connecting their incident vertices such that no edge passes through any other vertex. A drawing is \emph{simple} if any pair of edges intersect at most once - either in a common vertex or a proper \emph{crossing} in the relative interior of the edges. All drawings considered in this paper are simple and the term simple is mostly omitted. A drawing is \emph{plane} if it does not contain any crossing.

For a fixed integer $n$ let $D$ be a simple drawing of the complete graph $K_n$ and let $\mathcal {T}_D$ be the set of all drawings of plane spanning trees which are subdrawings of $D$. Note that $\mathcal {T}_D$ is non-empty, as it contains at least the $n$ stars in $D$ (where a \emph{star} contains all edges incident to a single vertex). Unless explicitly stated otherwise, the word \emph{tree} always refers to a plane spanning tree in $\mathcal {T}_D$, where the drawing $D$ is either clear from the context or the statement holds for any simple drawing of~$K_n$. Two (plane) subdrawings $H$ and $H'$ of a simple drawing $D$ are said to be \emph{compatible} if the union of $H$ and $H'$ is still plane.% subdrawing of $D$.

Let $\mathcal{F}(\mathcal {T}_D)$ be the (abstract) graph that has a vertex for each plane spanning tree in $\mathcal {T}_D$ and two vertices are adjacent 
if and only if the corresponding trees are compatible. We call $\mathcal{F}(\mathcal {T}_D)$ the \emph{compatibility graph} of $\mathcal {T}_D$. In this paper, we study properties of $\mathcal{F}(\mathcal {T}_D)$, focusing primarily on connectivity aspects:

\begin{ques}\label{ques:main}
	Let $n$ be an integer. Is the compatibility graph $\mathcal F(\mathcal {T}_D)$ connected for any simple drawing $D$ of the complete graph $K_n$?
\end{ques}

Note that the notion of compatibility is closely related to the notion of edge flips: An \emph{edge flip} in a plane spanning tree is the operation of removing an edge and replacing it with a new edge such that the resulting graph is again a plane spanning tree. In our setting, we further require this pair of edges to be non-crossing. In fact, one can simulate transformations via compatible trees in terms of crossing free edge flips:
for two compatible trees $T_1, T_2$, successively add edges from $T_2$ to $T_1$, 
while removing an edge that is not in $T_2$ from the resulting cycle.

\begin{figure}[t]
\centering
\includegraphics[page=1, scale=1]{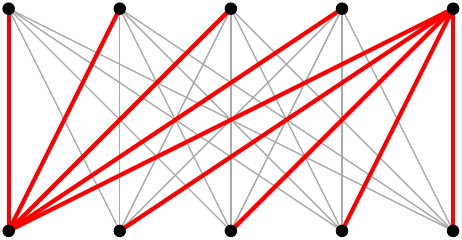}
\caption{A simple drawing of the complete bipartite graph with a tree (drawn in red, bold edges) that is an isolated vertex in the corresponding compatibility graph.}
\label{fig:bipartite}
\end{figure}

We observe that the compatibility graph of simple drawings 
that are not of the complete graph might not be connected even if the 
graph is dense. For example, Figure~\ref{fig:bipartite} shows a simple drawing of the complete bipartite graph containing a plane tree that crosses all edges 
of the graph not belonging to the tree. Hence, this tree is an isolated 
vertex in the correponding compatibility graph.

\paragraph*{\bf Related work.} The problem of transforming elements within a class of objects (e.g. plane spanning trees or matchings) into each other via a certain operation (e.g. edge flips or compatibility) has been studied extensively in a huge variety of contexts. Considering edge flips, some of the earliest results have been obtained on triangulations: Wagner~\cite{Wagner1936} showed connectivity of the corresponding flip-graph in the combinatorial setting and Lawson~\cite{lawson1972transforming} in the straight-line setting. For more details we refer the reader to the survey of Bose and Hurtado \cite{survey}.

Considering the notion of compatibility, most of the work has been done in the straight-line setting, e.g., in the context of perfect matchings with \cite{abhpv-ltdbm2018,abls2015} or without \cite{abdgh09,aght-cmgg-11} vertex coloring, or for edge-disjoint compatibility \cite{aam-dcgnc-15,ist-dcgm-13}. Aichholzer et al.~\cite{AICHHOLZER200619} showed, in the straight-line setting, that the compatibility graph of plane spanning trees is connected with diameter $O(\log k)$, where $k$ denotes the number of convex layers of the point set. Buchin et al.~\cite{buchin2009transforming} provided a corresponding worst case lower bound of $\Omega(\log n / \log \log n)$. 

It is natural to extend this question to simple drawings, which however are inherently difficult to handle (even the existence of certain plane substructures is still unresolved in simple drawings; see e.g. \cite{r-gdcg-88}). On the positive side, Garc\'{i}a, Pilz and Tejel~\cite{garcia2021plane} proved that any maximal plane subgraph is 2-connected, which guarantees for any plane spanning tree the existence of a compatible plane spanning tree. In this paper, we aim to shed some light on this wide open topic of compatibility graphs of trees in simple drawings.

\begin{figure}[t]
\centering
\includegraphics[page=1, scale=1]{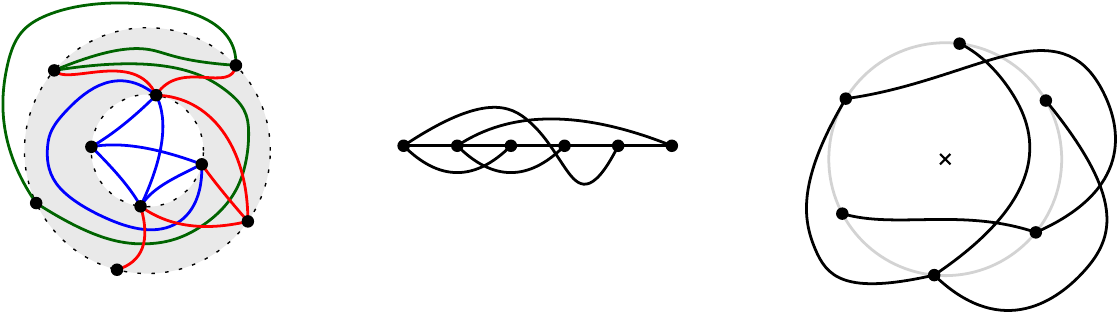}
\caption{\emph{Left to right:} cylindrical, monotone, strongly c-monotone drawing.}
\label{fig:basic_definitions}
\end{figure}

\paragraph*{\bf Contribution.} We approach Question~\ref{ques:main} from two directions, proving a positive answer for special classes of drawings (namely, cylindrical, monotone, and strongly c-monotone drawings) and for special classes of spanning trees (namely stars, double stars, and twin stars). We postpone the precise definitions of these classes of drawings and graphs to the later sections, however, \Cref{fig:basic_definitions} gives an illustration of these notions.

\begin{theorem}\label{thm:main_special_drawings}
	Let $D$ be a cylindrical, monotone, or strongly c-monotone drawing of the complete graph $K_n$. Then, the compatibility graph $\mathcal F(\mathcal {T}_D)$ is connected.
\end{theorem}

\begin{restatable}{theorem}{thmMainSpecialTrees}\label{thm:main_special_trees}
Let $D$ be a simple drawing of the complete graph $K_n$ and let $\mathcal {T}_D^*$ be the set of all plane spanning stars, double stars, and twin stars on $D$. Then, the compatibility graph $\mathcal F(\mathcal {T}_D^*)$ is connected.
\end{restatable}

\Cref{sec:specialdrawings} is devoted to the proof of \Cref{thm:main_special_drawings}, while \Cref{sec:special_trees} is dedicated to the proof of \Cref{thm:main_special_trees}. All results that are marked by a (clickable) $\textcolor{red}{(\star)}$ have a full proof in the appendix.

\section{Special simple drawings of $K_n$}\label{sec:specialdrawings}

In this section we prove connectedness of the compatibility graph for certain classes of drawings. Clearly, for any drawing of $K_n$ that admits a plane spanning tree which is not crossed by any edge of $D$, the compatibility graph is connected with diameter at most 2. This is, for example, the case for 2-page book drawings, where the vertices are placed along a line and each edge lies entirely in one of the two open halfplanes defined by this line.

\subsection{Cylindrical drawings}

Following the definition of Schaefer \cite{Schaefer2013},
in a \emph{cylindrical drawing} of a graph the vertices are placed along two concentric circles, the \emph{inner} and \emph{outer} circle, and no edge is allowed to cross these circles.

\begin{restatable}{lemma}{lemCylindrical}\label{lem:cylindrical}
\hyperref[lem:cylindrical:proof]{($\star$)}
	Let $D$ be a cylindrical drawing of $K_n$.  
	Then $\mathcal F(\mathcal {T}_D)$ is connected with diameter at most 4. 
\end{restatable}

\subsection{Monotone drawings}\label{sec:monotone}

A simple drawing in which no two vertices have the same $x$-coordinate and every edge is drawn as an $x$-monotone curve is called \emph{monotone drawing}. Let $v_1,v_2, \ldots , v_{n}$ denote the sequence of vertices in increasing $x$-order. W.l.o.g. assume that these vertices are on the $x$-axis. Then, the plane spanning path $\mathcal{S} = v_1,v_2, \ldots , v_{n}$ is called \emph{spine} path. 
An edge that intersects the spine path is called \emph{twiggly} edge. 

We define a relation on the twiggly edges of $D$ as follows: 
	for two twiggly edges $e,f$ we have $e \succ f$  
	if they are non-intersecting and admit a vertical line intersecting the relative interiors of both edges that intersects $e$ at a larger $y$-coordinate than $f$. All other pairs of twiggly edges are incomparable.
For a set~$E$ of pairwise non-intersecting twiggly edges, an edge $e \in E$ is \emph{maximal} if there is no other edge $f\in E$ s.t. $f \succ e$. Note that this relation is acyclic, i.e., there are no twiggly edges $e_1, \ldots, e_k$ such that $e_1 \succ e_2 \succ \ldots \succ e_k \succ e_1$. And hence, any non-empty set of twiggly edges admits a maximal element.

\begin{restatable}{lemma}{lemMonotone}\label{lem:monotone}
	\hyperref[lem:monotone:proof]{($\star$)}
For any monotone drawing $D$ of $K_n$, the compatibility graph $\mathcal F(\mathcal {T}_D)$ is connected with diameter $O(n)$.
\end{restatable}

\begin{figure}[t]
	\centering\includegraphics[page=1, scale=0.4]{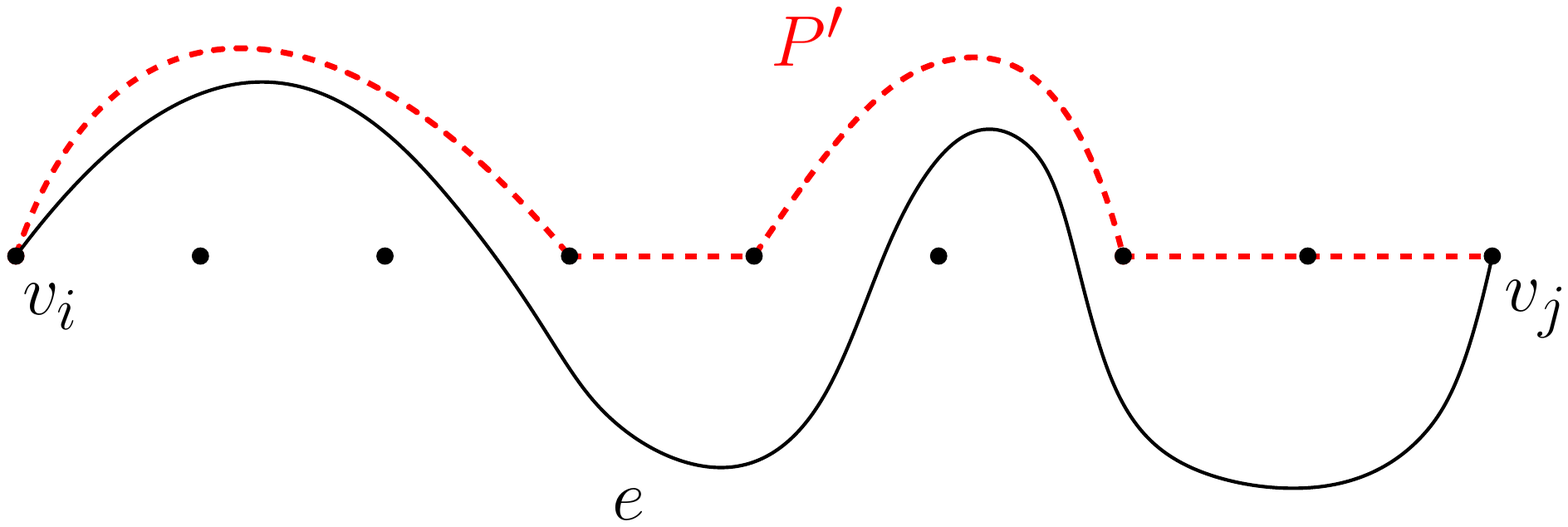}
    \caption{The (maximal) twiggly edge $e=v_iv_j$ divides the vertices between $v_i$ and $v_j$ into two groups -- above and below. 
	The path $P'$ is formed by joining the consecutive vertices lying above $e$ including the vertices of $e$. 
	}  \label{fig:monotonecase}
\end{figure}

\begin{proof}[Sketch]
We show that any plane spanning tree $T$ in $D$ can be transformed to the spine path $\mathcal{S}$. If $T$ does not contain any twiggly edge, clearly it is compatible to $\mathcal{S}$. 
Otherwise, we proceed as follows. Corresponding to a maximal twiggly edge $e$ of $T$, we find a path $P'$ connecting the vertices of $e$ (see \Cref{fig:monotonecase}). We can show that $P'$ is compatible to $T$ and lies strictly above $e$. Thus, we can add $P'$ to $T$, which creates at least one cycle in $T$. Removing appropriate edges including $e$, we get a compatible tree with at least one twiggly edge less and repeating this process, we will eventually reach the spine path $\mathcal{S}$. 
\end{proof}

\subsection{Strongly c-monotone drawings}

A curve is called \emph{c-monotone} (w.r.t. a point $x$) if every ray emanating from $x$ intersects the curve at most once.
A simple drawing is \emph{c-monotone}, if all vertices are drawn along a circle and every edge is a c-monotone curve w.r.t. the center of the circle. A c-monotone drawing is \emph{strongly} c-monotone if for any pair of edges $e,e'$ there is a ray (rooted at the circle center) that neither intersects $e$ nor $e'$.

In a (strongly) c-monotone drawing, we label the vertices $v_1,v_2, \ldots, v_n$ in cyclic order and denote the center of the circle by $c$.
In the following, we often consider edges and their intersections with rays rooted at $c$; unless stated otherwise, any ray is rooted at $c$ and edges are intersected in their relative interiors.

An edge $e$ connecting two consecutive vertices $v_i$, $v_{i+1}$ is called \emph{cycle edge} and if $e$ is drawn along the \enquote{shorter} side of the circle it is called \emph{spine edge} (that is, no ray formed by the center and any vertex intersects~$e$). All spine edges form the \emph{spine} and any path consisting entirely of spine edges is called \emph{spine path}.

\begin{restatable}{lemma}{lemSpineEdges}\label{lem:spine_edges}
\hyperref[lem:spine_edges:proof]{($\star$)}
Any strongly c-monotone drawing $D$ of $K_n$  either has all cycle edges as spine edges or is isomorphic to a monotone drawing.
\end{restatable}

Again, we define \emph{twiggly} edges to be those that intersect a spine edge. 
A crucial difference to the monotone setting is that an analogue to the relation '$\succ$' (adjusted with respect to the intersection with rays emanating from $c$) may now be cyclic and hence, we cannot guarantee the existence of a \emph{maximal} twiggly edge anymore. We therefore need a different approach.

\begin{figure}[t]
	\centering\includegraphics[page=1, scale=0.6]{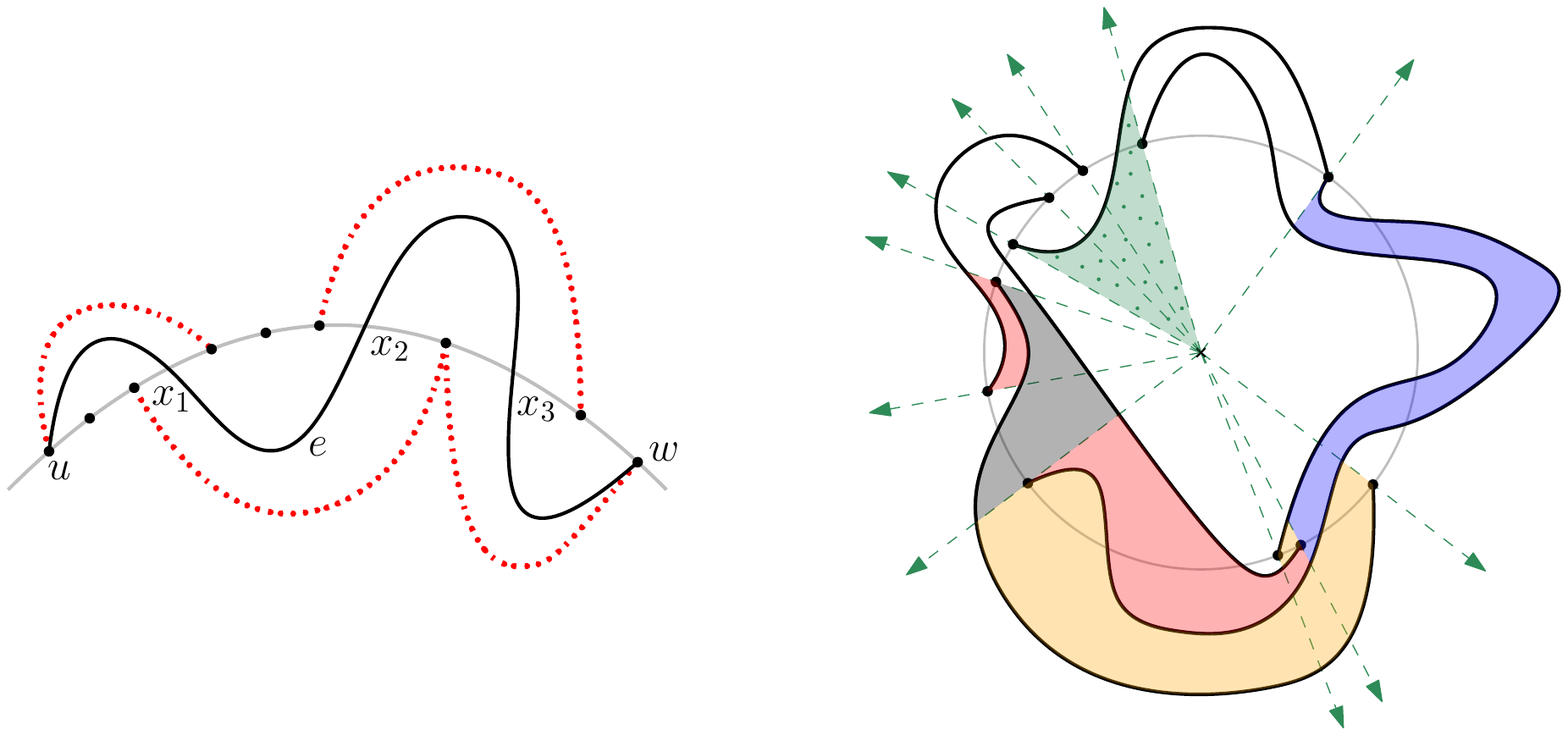}
	\caption{\emph{Left:} The red dotted edges are bumpy edges of the twiggly edge~$e$. \emph{Right:} A set of twiggly edges and some corridors; the dotted green is an inner corridor.}  \label{fig:c_monotone_notions}
\end{figure}

For a twiggly edge $e = uw$, let $x_1, \ldots, x_k$ be its crossings with the spine (note that these are not vertices of $K_n$) and assume the labeling to be in such a way that $u, x_1, \ldots, x_k, w$ appear in clockwise order. For $i \in \{1,\ldots, k\}$ denote the vertex (of $K_n$) in clockwise order before~$x_i$ by~$x_i^-$ and the one after by $x_i^+$. Furthermore, set $u = x_0^-$ and $w = x_{k+1}^+$. Then, for $i \in \{0, \ldots, k\}$, we call the edges $x_i^-x_{i+1}^+$ \emph{bumpy} edges (see~\Cref{fig:c_monotone_notions} (left)). Note that bumpy edges do not intersect the spine and for any twiggly edge there are at least two bumpy edges.

Clearly, we can identify any ray $r$ with an angle $\theta$, the angle it forms with the vertical ray (upwards). Two edges $e,f$ are called \emph{neighbours on an interval $[\theta_1, \theta_2]$}, if for any ray $r \in [\theta_1, \theta_2]$ the intersections of $e$ and $f$ with $r$ appear consecutively on $r$. A \emph{corridor} is a maximally connected region bounded by two neighbouring edges (along a maximal interval). Again, we identify corridors by an interval $[\theta_1, \theta_2]$ and usually we speak of corridors defined by the edges of a plane spanning tree. The \emph{twiggly depth} (with respect to a plane spanning tree~$T$) of a ray $r$ is the number of twiggly edges (of $T$) that $r$ intersects.

We extend our definition of neighbours (along an interval) also to the very inside and very outside by inserting a dummy edge at the circle center and one at infinity. 
More precisely, an edge $e$ is the neighbor of the circle center $c$ along an interval $[\theta_1, \theta_2]$ if for any ray $r \in [\theta_1, \theta_2]$ the intersection of $r$ and $e$ is closest to $c$ (and furthest in the case of being a neighbor of infinity). We call the corresponding corridors \emph{inner/outer} corridors. Note that the set of all corridors partitions the plane. See \Cref{fig:c_monotone_notions} (right) for an illustration.

We further remark that for any plane spanning tree $T$, any corridor $C=~[\theta_1, \theta_2]$ (of edges of~$T$) begins and ends at a vertex, i.e., the rays at $\theta_1$ and $\theta_2$ hit a vertex.

\begin{restatable}{lemma}{lemCorridor}\label{lem:corridor}
\hyperref[lem:corridor:proof]{($\star$)}
For any plane spanning tree $T$ of a strongly $c$-monotone drawing~$D$ and any corridor $C$ of $T$ with start and end vertex $s$ and $t$, there is a path~$P$ in $D$ from $s$ to $t$ staying entirely in $C$, that does not intersect $T$.
Furthermore, if $C$ is an inner or outer corridor, $P$ does not use any twiggly edge.
\end{restatable}

\begin{lemma}\label{lem:c_monotone}
For any strongly c-monotone drawing of $K_n$, the compatibility graph $\mathcal F(\mathcal {T}_D)$ is connected with diameter $O(n)$.
 
\end{lemma}

\begin{proof}
Let $D$ be a strongly c-monotone drawing of $K_n$ and let $T$ be a plane spanning tree. We show that $T$ can be compatibly transformed to a spine path (by iteratively decreasing its twiggly depth). By \Cref{lem:spine_edges}, we may assume that all $n$ spine edges are present in $D$. Again, if there is no twiggly edge in $T$, then $T$ is compatible with the spine.

Let $E_{twig}$ be the set of twiggly edges of $T$ and construct the set $\mathcal{C}$ of all corridors. Next, for any corridor $C \in \mathcal{C}$ with start and end vertex $s$ and $t$, we add the path $P_C$ as guaranteed by \Cref{lem:corridor} to $T$. 

Clearly, we do not disconnect $T$ when removing $E_{twig}$ now. Indeed, let $e = uw \in E_{twig}$, then the collection of corridor paths below (and also above) $e$ connects $u$ and $w$. 
	So we remove $E_{twig}$ and potentially some further edges  
	until $T$ forms a spanning tree again (which by \Cref{lem:corridor} is also plane). 
	Furthermore, any ray $r$ that intersects $x$ previous twiggly edges (i.e., $E_{twig}$) intersects $x+1$ corridors, two of which are either an inner or outer corridor. By \Cref{lem:corridor} and the properties of $c$-monotone curves, $r$ intersects at most $x-1$ (new) twiggly edges. Hence, the twiggly depth of any ray decreased by at least one and we recursively continue this process until all rays have twiggly depth 0, in which case $T$ is compatible to a spine path. As we have twiggly depth at most $n-1$ in the beginning, 
$\mathcal F (\mathcal{T}_D)$ has diameter $O(n)$.
\end{proof}

\Cref{thm:main_special_drawings} now follows from \Cref{lem:cylindrical}, \Cref{lem:monotone}, and \Cref{lem:c_monotone}.

\section{Special Plane Spanning Trees}\label{sec:special_trees}

In this section, we are not restricting our drawing anymore, i.e., $D$ will be a simple drawing of $K_n$ throughout this section. Instead we focus on special classes of spanning trees and show that the subgraph $\mathcal F (\mathcal {T}_D^*)$ of $\mathcal F(\mathcal {T}_D)$ induced by the set of vertices corresponding to stars, double stars, and twin stars is connected.

A plane spanning tree with a fixed path $P$ of length $k$ such that all other vertices are incident to either the start or end vertex of $P$ is called a \emph{$k$-star}. 
A $0$-star (i.e., $P$ consists of a single vertex) is called \emph{star}. A $1$-star is called \emph{double star} and a $2$-star is called \emph{twin star}. 

The following relation, introduced in~\cite{apsvw-sssdk-19}, will be very useful: Given a simple drawing of $K_n$ with vertex set $V$ and two vertices $g\neq r \in V$, for any two vertices $v_i, v_j \in V \backslash\{g,r\}$, we define $v_i \rightarrow_{gr} v_j$ if and only if the edge $v_ir$ crosses $v_jg$. In~\cite{apsvw-sssdk-19} it is shown that this relation is asymmetric and acylic.

We start by showing that stars can always be transformed into each other via a sequence of crossing free edge flips.

\begin{restatable}{lemma}{lemStarmain}\label{thm:stars}
	
	Any two stars in $D$ have distance $O(n)$ in $\mathcal F (\mathcal {T}_D^*)$.
\end{restatable}

\begin{figure}[t]
	\centering\includegraphics[page=1, scale=0.38]{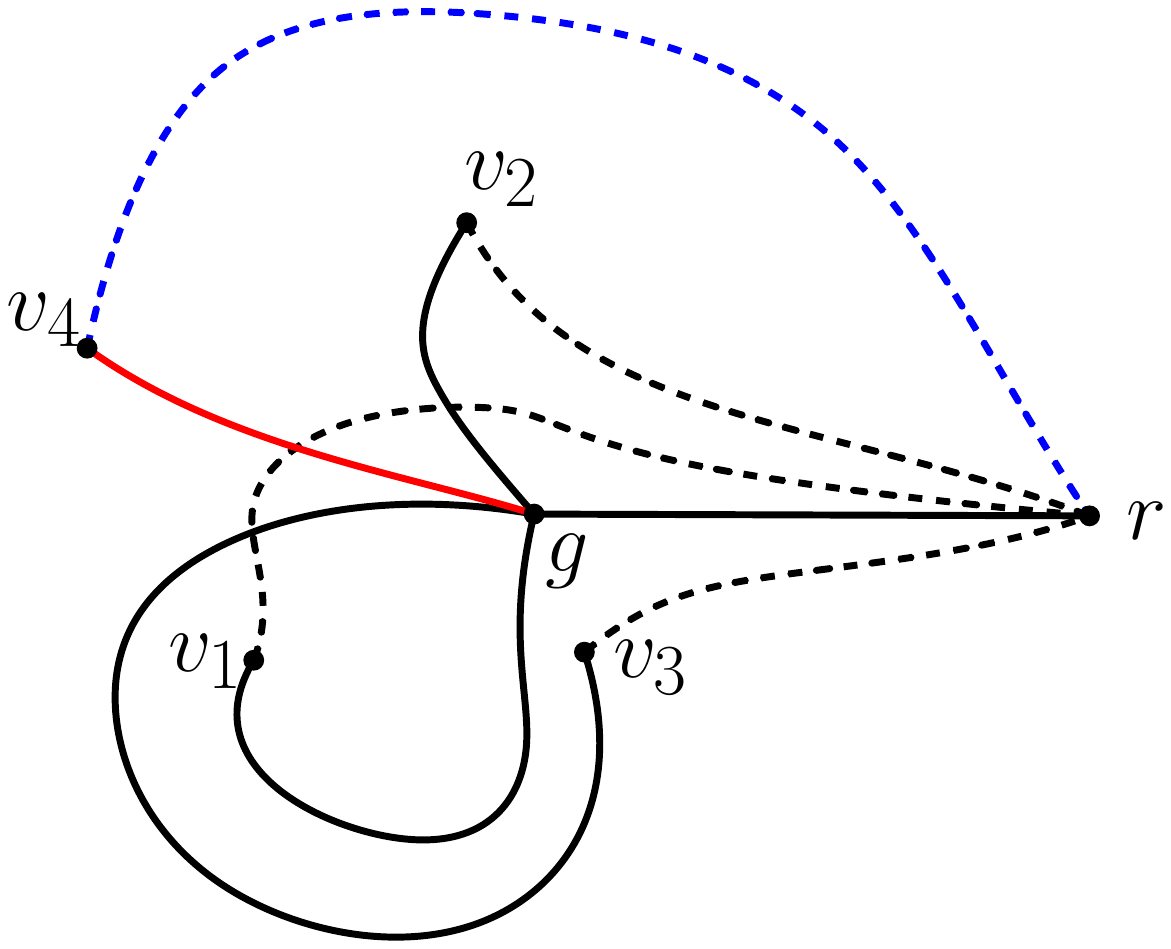}
	\caption{Proof of \Cref{thm:stars}: The solid edges represent a star in $g$, while the dotted edges form a star in $r$. The vertices are labeled conforming to the relation $\rightarrow_{gr}$. In order to transform the star in $g$ to the star in $r$, the first step is adding the dotted blue edge $v_4r$ and deleting the red edge $v_4g$.
	}  \label{fig:star}
\end{figure}

\begin{proof}

	Given a star $T$ in $g$ (i.e., $g$ is incident to all other vertices of $T$), we can transform it into a star $H$ in $r$ via a sequence of crossing free edge flips, such that in every step, the graph is a double star with fixed path $r$,$g$, in the following way. We label the vertices in $V \backslash\{g,r\}$ such that $v_i \rightarrow_{gr} v_j$ implies $i < j$ (see \Cref{fig:star}). We iteratively replace an edge $gv_i$ by $rv_i$ starting from $i = n-2$ and continuing in decreasing order. Clearly, all intermediate trees are double stars (with fixed path $r$,$g$) and hence, it remains to argue that the flips are compatible, i.e., for $i = n-2, \ldots, 1$ the edge $gv_i$ does not cross any edge of the current $T$.
By construction, in any step $i$, $T$ contains edges of the form (a) $rv_j$ for $j > i$ and (b) $gv_k$ for $k < i$. The edge $gv_i$ cannot cross edges in (a) by the definition of the relation $\rightarrow_{gr}$ and also not those in (b) due to the properties of simple drawings. As we need at most $n-2$ steps for the transformation, any two stars have distance $O(n)$ in $\mathcal F (\mathcal{T}_D^*)$.
\end{proof}

\Cref{thm:main_special_trees} then follows from \Cref{thm:stars} in combination with the following two lemmata.

\begin{restatable}{lemma}{lemStarDoublestar}\label{lem:star_doublestar}
\hyperref[lem:star_doublestar:proof]{($\star$)}
Any double star in $D$ has distance $O(n)$ to any star in $\mathcal F (\mathcal {T}_D^*)$.
\end{restatable}

\begin{restatable}{lemma}{lemStarTwinstar}\label{lem:star_twinstar}
\hyperref[lem:star_twinstar:proof]{($\star$)}
Any twin star in $D$ has distance $O(n)$ to any star in $\mathcal F (\mathcal {T}_D^*)$.
\end{restatable}

\bibliographystyle{splncs04}
\bibliography{spanning_trees_in_simple_drawings}

\appendix

\newpage

\section{Missing proofs of Section~\ref{sec:specialdrawings}}\label{app:sec:specialdrawings}

We need a few more definitions and observations for the proof of \Cref{lem:cylindrical}.
We denote the vertices on the inner/outer circle of a cylindrical drawing by \emph{inner/outer vertices} and the edges connecting any two inner/outer vertices by \emph{inner/outer edges}. The remaining edges are called \emph{side edges}. Furthermore, edges that join consecutive vertices of either of the two circles are called \emph{cycle edges}.

\begin{remark}\label{rem:cycle_edges}
Cycle edges may only be crossed by side edges and either of the two circles contains at most one such cycle edge that is crossed.
\end{remark}

In other words, along both circles there is a Hamiltonian path of cycle edges that is not crossed by any other edge.

Further note that any plane spanning tree in a cylindrical drawing contains at least one side edge.

\lemCylindrical*

\begin{proof}\label{lem:cylindrical:proof}

\emph{(Disclaimer: we moved a few definitions and remarks, that are relevant for this proof, to the preceeding paragraphs.)}
Let $T_1,T_2 \in \mathcal {T}_D$ and $e_1 = v_sv_r$ be a side edge of $T_1$ as well as $e_2 = v_s'v_r'$ be a side edge of $T_2$. Furthermore, as guranteed by \Cref{rem:cycle_edges}, let $S_1 \in \mathcal {T}_D$ be a plane spanning tree consisting of the two Hamiltonian paths of cycle edges and $e_1$ (and similarly $S_2 \in \mathcal {T}_D$ consists of the cycle paths and $e_2$). Clearly, $T_1$ and $S_1$ are compatible as well as $T_2$ and $S_2$. If $S_1$ and $S_2$ are compatible, we are done. Otherwise, $e_1$ and $e_2$ are crossing in $D$. As the induced 4-tuple $\{v_s,v_s',v_r,v_r'\}$ has at most one crossing, the two side edges $v_sv_r'$ and $v_s'v_r$ cannot cross $e_1$ and $e_2$. Now consider one of these edges as $e'$ and construct a plane spanning tree $S_3$, using only $e'$ and the uncrossed cycle paths. Clearly, $S_3$ is compatible to $S_1$ as well as $S_2$ and hence, $T_1$,$S_1$,$S_3$,$S_2$,$T_2$ is a path of length 4 in $\mathcal F(\mathcal {T}_D)$.  
\end{proof}	

\lemMonotone*

\begin{proof}\label{lem:monotone:proof}
	We prove the connectedness by showing that any plane spanning tree in $D$ can be transformed (via a sequence of compatible trees) to the spine path $\mathcal{S}$. So, let $T \in \mathcal {T}_D$ be a plane spanning tree. If $T$ does not contain any twiggly edge, then it is compatible to $\mathcal{S}$. 
	Otherwise, let $e = v_i v_j$ (with $i<j$) be a maximal twiggly edge of $T$. Note that all twiggly edges in $T$ are non-intersecting, since $T$ is plane. Define
	\[
	V_e^{\uparrow} = \{ v_k : i < k < j, v_k \text{ is above } e\}.
	\]

	Consider the path $P' = v_i, V_e^{\uparrow}, v_j$ that starts at $v_i$, ends at $v_j$ and inbetween uses only vertices of $V_e^{\uparrow}$ (see \Cref{fig:monotonecase}). Note that the vertices in $P'$ are in increasing $x$-order. Clearly, no edge in $P'$ can intersect the twiggly edge $e$, since it would have to intersect it twice. Hence, $P'$ lies (strictly) above $e$. Next we show that $P'$ cannot intersect $T$. To this end, let $f \in T$ and assume, for the sake of contradiction, that it intersects an edge $v_xv_y$ of $P'$. Clearly, $f$ can neither be incident to $v_x$ nor $v_y$. However, by construction of $P'$, in order to reach its incident vertices, $f$ has to either (i) intersect $v_xv_y$ twice, or (ii) intersect $e$, or (iii) intersect a spine edge. (i) yields a contradiction to the properties of simple drawings, (ii) yields a contradiction to the planarity of $T$, and (iii) yields a contradiction to the maximality of the twiggly edge $e$.
	
	Hence, we can add $P'$ to $T$ without introducing any crossings. Note that some edges of $P'$ may already be present in $T$, however, since this insertion creates a cycle, at least one is not. In order to reach a compatible tree, we remove $e$ from this cycle (and potentially more edges until reaching a plane spanning tree again). We created a compatible tree with at least one twiggly edge less and repeating this process, we will eventually reach the spine path $\mathcal{S}$. As we have at most $n-1$ twiggly edges, $\mathcal F (\mathcal{T}_D)$ has diameter $O(n)$.
\end{proof}

\lemSpineEdges*

\begin{proof}\label{lem:spine_edges:proof}
Assume there is an edge $e = v_iv_{i+1}$ of consecutive vertices which does not form a spine edge, i.e., $e$ intersects any ray which is not in the wedge $v_icv_{i+1}$. The c-monotonicity of edges implies that any edge intersecting the wedge $v_icv_{i+1}$ must intersect any ray in this wedge. This, however, would be in contradiction to the property of strongly c-monotone drawings and hence, the wedge $v_icv_{i+1}$ is not intersected by any edge. Clearly, this implies the drawing to be (strongly) isomorphic to a monotone drawing: any ray $r$ in the wedge $v_icv_{i+1}$ is uncrossed by all edges and hence, we can cut the spine at the intersection with $r$ and stretch it to a monotone drawing.
\end{proof}

\lemCorridor*

\begin{figure}[t]
	\centering\includegraphics[page=1, scale=0.5]{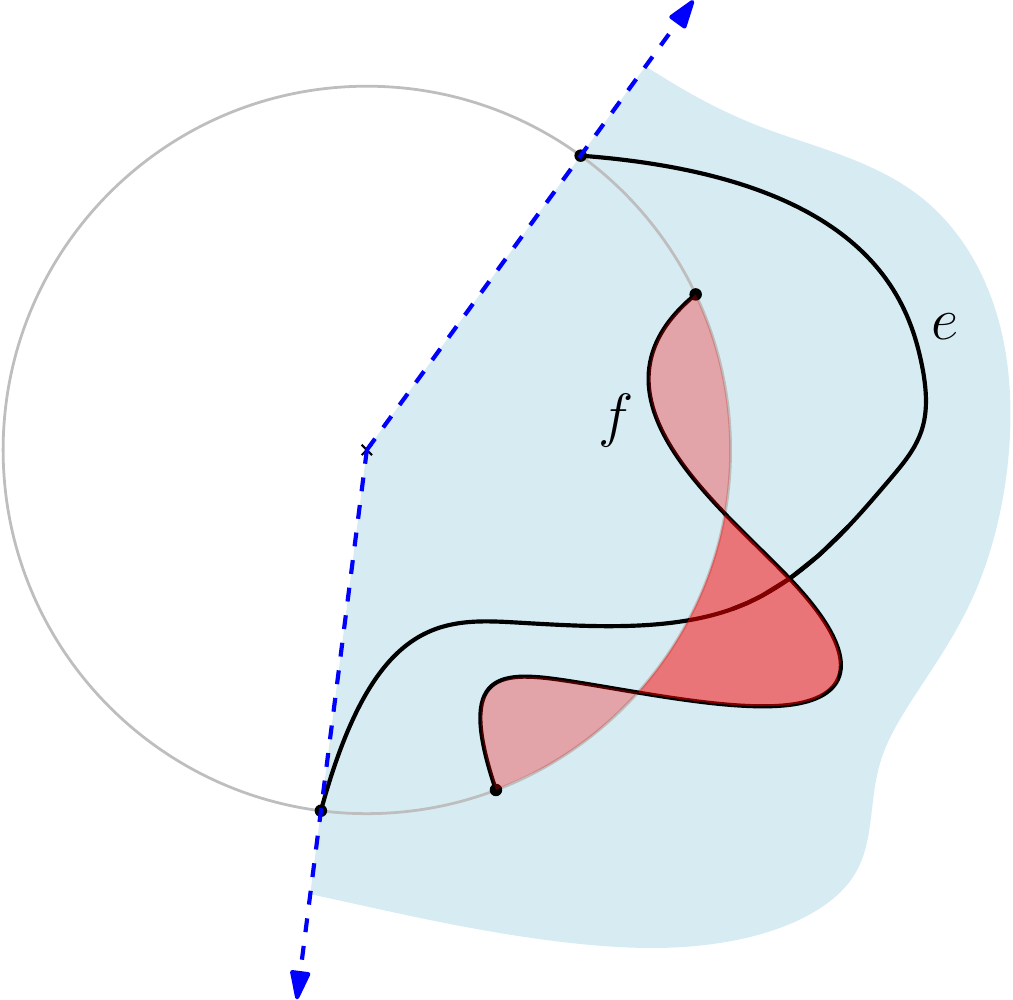}
	\caption{The faces formed by $f$ and the spine are red. If $e$ intersects $f$ it must leave the corresponding face by intersecting the spine.}
	\label{fig:surrounding}
\end{figure}

\begin{proof}\label{lem:corridor:proof}
Start the path $P$ at $s$ and always connect to the next vertex (along the circle) that is in the corridor $C$. Due to strong c-monotonicity, the corresponding edge must run through $C$. Note that $P$ may consist of bumpy, spine, and twiggly edges. However, $P$ may only contain a twiggly edge if it is forced to by \emph{both} edges bounding the corridor (e.g., the blue corridor in \Cref{fig:c_monotone_notions} (right)).

Moreover, since we always connect consecutive vertices along the corridor, any edge $e$, intersecting an edge $f$ of $P$, must \emph{surround}\footnote{We say that an edge $e = uw$ \emph{surrounds} an edge $f$ (all four vertices are distinct), if the wedge $ucw$ (containing $e$) also contains $f$; see \Cref{fig:surrounding}.} $f$.
We claim that $e$ must also be a twiggly edge, which can be seen as follows. The twiggly edge $f$ together with the spine forms a collection of faces which are only bounded by $f$ and the spine (see \Cref{fig:surrounding}). Since $e$ surrounds $f$, $e$ cannot have a vertex within such a face. Hence, if $e$ intersects $f$, $e$ must also intersect the spine and thus, is a twiggly edge. Therefore, $e$ cannot belong to $T$ (since it is in the corridor). In particular, $P$ does not intersect~$T$.

Lastly observe that the paths in an inner/outer corridor do not use any twiggly edge, since any inner/outer corridor is bounded by a dummy edge which does not cross the spine.
\end{proof}

\section{Missing proofs of Section~\ref{sec:special_trees}}\label{app:sec:special_trees}

\lemStarDoublestar*

\begin{proof}\label{lem:star_doublestar:proof}
Let $T$ be a double star with fixed path $g$,$r$ and $T'$ be a star in $r'$. $T$ can be transformed to a star in $r$ (or $g$) completely analogous as in the proof of \Cref{thm:stars}. The only difference is that there are additional edges attached to $r$, i.e., edges of type (a) in \Cref{thm:stars}, which do not interfere. Next, using \Cref{thm:stars} again, we transform this star in $r$ to $T'$. This also implies the distance between a double star and a star in $\mathcal F (\mathcal{T}_D^*)$ to be $O(n)$.
\end{proof}

\lemStarTwinstar*

\begin{proof}\label{lem:star_twinstar:proof}
Let $T$ be a twin star with fixed path $g,s,r$ and $T'$ be a star in $r'$. Note that all edges in $T$ are incident to $g$ or $r$ and hence, we can add the edge $gr$ and remove $gs$ or $rs$ in order to obtain a double star. Using \Cref{lem:star_doublestar} this double star can be transformed to $T'$. Thus, a twin star can be tranformed to a star in $O(n)$ steps.
\end{proof}

\end{document}